\documentclass[12pt,twoside]{article}
\usepackage[all]{xy}
\usepackage{a4,amsmath,amssymb,amsfonts,amscd,mathrsfs}
\reversemarginpar
\newcommand{\kkk}[1]{}

\newcounter{Abschnitt}[section]
\newcommand{\neu}[1]{\protect\refstepcounter{Abschnitt}\protect
   \label{#1}\vspace{1ex}
   {\bf \protect\arabic{section}.\protect\arabic{Abschnitt}}
                     \kkk{#1}}

\newcommand{\rk}{{\rm rk}}

\newcommand{\Grass}{{\rm Grass}}

\newcommand{\Tor}{{\rm Tor}}
\newcommand{\Spec}{{\rm Spec}}

\newcommand{\Sym}{{\rm Sym}}

\newcommand{\Ccal}{{\cal C}}

\newcommand{\Jcal}{{\cal J}}

\newcommand{\Mcal}{{\cal M}}
\newcommand{\Ocal}{{\cal O}}

\newcommand{\cdop}{{\mathbb C}}

\newcommand{\pdop}{{\mathbb P}}

\newcommand{\dual}{^\lor}

\newcommand{\rarpa}[1]{\stackrel{#1}{\rightarrow}}

\newcommand{\proof}{{\bf Proof: }}

\newcommand{\qed}{{ \hfill $\square$}}

\newcommand{\tip}{{\vspace{0.5em} }}

\author{Georg Hein}
\title{Residual intersection theory with reducible schemes}
\begin{document}
\maketitle
\begin{abstract}
We develop a formula (Theorem \ref{Smain}.\ref{main})
which allows to compute top Chern classes
of vector bundles on the vanishing locus $V(s)$ of a section of this bundle.
This formula particularly applies in the case when
$V(s)$ is the union of locally complete intersections
giving the individual contribution of each component and their
mutual intersections.
We conclude with applications to the enumeration of rational curves
in complete intersections in projective space.
\end{abstract}

\section{Introduction}
We consider a vector bundle $E$ of rank $n$ on a scheme $X$ of dimension $n$.
If $E$ is globally generated, then a general section $s$ of $E$ vanishes in
$\int_X c_n(E)$ points.
This gives an easy way to compute the top Chern number of $E$.

However, if $E$ is not globally generated or if we have just one section of $E$
with a higher dimensional vanishing locus $V(s)$,
then we can compute this Chern number on this vanishing locus.
This was carried out in \cite{Ful}, and has a particularly nice form if $V(s)$ is a locally complete intersection (see Proposition 9.1.1 in \cite{Ful}).

We concentrate on the case where $V(s) = \cup_{i=1}^M Z_i$ is a union of locally complete intersections, which intersect mutually in a nice way (see \ref{Slocal}.\ref{local1} for more details).
Our main result is Theorem \ref{Smain}.\ref{main} where we give the contribution of each $Z_{i_1} \cap Z_{i_2} \cap \ldots \cap Z_{i_k}$ to the top Chern number of $E$.
In section \ref{Scounting} we outline how our result may be used to count rational curves on complete intersections.
We reduce the problem to computations on moduli spaces of rational curves in $\pdop^n$ and show that the corresponding vanishing cycle satisfies our assumption \ref{Slocal}.\ref{local1}.
Finally, we give two examples how our theorem might be used to compute some Chern numbers.
In the first example we compute the number of lines on a smooth cubic surface in $\pdop^3$.
Here a point, which is the intersection of two irreducible components,
gives a contribution of $-6$,
whereas most people would expect its contribution to be $\pm 1$.
In the second example we compute the number of lines on a general quintic threefold in $\pdop^4$ by linearizing the problem.
To compute these numbers Kontsevich's idea of using a $\cdop^*$-action and  Bott's fix point formula gives a more efficient way to compute these numbers along smooth subvarieties (see \cite{Kon}).
However, when there is no $\cdop^*$-action our technique might be used as
an instrument to reduce the dimension in these cases.

\section{Notations and preliminaries}\label{Snotations}
\neu{Segre}{\bf Segre classes.}
In the sequel it will turn out that for our problem it is convenient to
work with Segre classes, even though Chern classes are favored by most
mathematicians. The $k$-th Segre class of a vector bundle $E$ of rank $r$ on $X$ is
a homomorphism $A_{\bullet +k }(X) \to A_{\bullet}(X)$ of the group of
rational cycle classes. It is defined by
$$s_k(E)(\alpha) := (-1)^k \pi_*(c_1(\Ocal_E(1))^{r-1+k} \cap \pi^*(\alpha)) $$
where $\pi:\pdop(E) \to X$ is the projective bundle associated to $E$,
and $\Ocal_E(1)$ the relative ample quotient linebundle on $\pdop(E)$.
This definition coincides with the standard definition given in \cite{Ful}.
The sign which does not appear in \cite{Ful} has to be introduced because we use $\pdop(E)$ instead of $\pdop(E\dual)$.

\neu{multinomial}{\bf Multinomial coefficients.}
If $d=d_1+d_2+\ldots d_k$, then multinomial coefficient
$\binom{d}{d_1,d_2,\ldots,d_k}$ is defined to be the integer
$\frac{d!}{d_1! \cdot d_2! \cdot \ldots d_k!}$.
It is the coefficient of the monomial
$x_1^{d_1}\cdot x_2^{d_2} \cdot \ldots \cdot x_k^{d_k}$ in $(x_1+x_2+\ldots + x_k)^d$.
These coefficients satisfy the basic equation 
$\binom{d}{d_1,d_2,\ldots,d_k} = \binom{d-1}{d_1-1,d_2,\ldots,d_k}
+ \binom{d-1}{d_1,d_2-1,\ldots,d_k} + \ldots + 
\binom{d-1}{d_1,d_2,\ldots,d_k-1}$.

\section{Reduction to the vanishing locus $V(s)$}\label{Sreduction}
\neu{reduction1}
Let $X$ be a smooth scheme of dimension $n$,
and $E$ be a vector bundle on $X$ of rank $n$.
We want to compute the integral $\int_X c(E)$, 
or the $n$-th Chern number of $E$.
However, we want to do all these computations
on the vanishing scheme $V(s)$ of a global
section $s \in H^0(E)$.
This allows the reduction to lower dimensional schemes.
The ideal of $V(s)$ is given by the image of the dual
morphism to $s: \Ocal_X \to E$.
We distinguish the following three cases:

{\bf Case 1.} $s=0$ which implies $V(s)=X$.\\
Thus, we can work on all of $X$ and everything is easy.
However, we are not able to reduce to a lower dimensional scheme.
Therefore, we assume from now on that $s$ is not the trivial section.

{\bf Case 2.} $V(s)=:D$ is a Cartier divisor.\\
We obtain an exact sequence $0 \to \Ocal(D) \to E \to E' \to 0$.
We have therefore
$$\int_X c(E) = c(E')[D] = \int_D c(E'|_D) =
\int_D c(E|_D).c(\Ocal_X(D)|_D)^{-1} \, .$$

{\bf Case 3.} The general case.\\
We consider the blow up $\pi:\tilde X \to X$ of $X$ in $V(s)$.
Let $D$ be the exceptional divisor.
$$\xymatrix{D \ar[r] \ar[d]^\pi & \tilde X \ar[d]^\pi\\
V(s) \ar[r] & X}$$
The section $s$ gives rise to a section $\tilde s$ of $\pi^*E$
with vanishing locus $D$.
Thus, we can apply the  previous result and find

\neu{reduction2}
{\bf Proposition.} {\em (cf. Proposition 9.1.1 in \cite{Ful})}
$$\int_Xc(E) = \int_{\tilde X} c(\pi^*E) = 
\int_D c(\pi^*E|_D).c(\Ocal_{\tilde X}(D)|_D)^{-1}
=\int_{V(s)}c(E|_{V(s)}).\pi_*(c(\Ocal_{\tilde X}(D)|_D)^{-1}) \,.  $$

\neu{reduction3}
If $V(s)$ is a locally complete intersection of codimension $d$,
then $D$ is a $\pdop^{d-1}$-bundle over $V(s)$.
To be precisely, we have $D=\pdop(\Jcal / \Jcal^2)$.
The normal bundle $\Ocal_{\tilde X}(D)$ of $D$ in $\tilde X$
restricted to $D$ is the linebundle $\Ocal(-1)$.
Thus, we find
$$\pi_* (c(\Ocal_{\tilde X}(D)|_D)^{-1}) =
\pi_* (\sum_{k\geq 0} c_1^k(\Ocal_D(1) )) = s_{-1}(\Jcal/\Jcal^2)
= s(N_XV(s)) = c(N_XV(s))^{-1}\, .$$

Thus, in this case we obtain

\neu{reduction4}
{\bf Proposition.}
{\em (cf. Proposition 9.1.1 in \cite{Ful})
If $V(s)$ is a locally complete intersection, then we have\\
$$\int_X c(E) = \int_{V(s)} c(E|_{V(s)}) \cdot c(N_XV(s))^{-1} \, .$$
}

\neu{exam1}
{\bf Example: The Euler characteristic.}
Let $X$ be a smooth complex variety then we have $\chi(X)= \int_X c(T_X)$.
Let us assume that there exists a global section $s$ of $T_X$ such
that $V(s)$ is a locally complete intersection.
Using the result of Proposition \ref{Sreduction}.\ref{reduction4} we obtain
$\chi(X) = \int_{V(s)} c(T_X|_{V(s)}) \cdot c(N_XV(s))^{-1}$.
However, the standard exact sequence of the normal bundle
$0 \to T_{V(s)} \to T_X|_{V(s)} \to N_XV(s) \to 0$ yields
$c(T_{V(s)})= c(T_X|_{V(s)}) \cdot c(N_XV(s))^{-1}$.
Thus, we obtain $\chi(X)=\chi(V(s))$.

\neu{exam2}{\bf The Euler characteristic of Grassmannians.}
Let $\Grass(n,k)$ be the Grassmannian scheme of rank $k$ quotients of $\cdop^n$.
On $\Grass(n,k)$ we have the tautological short exact sequence
$0 \to K \rarpa{\alpha} \cdop^n \otimes \Ocal_{\Grass(n,k)} \to Q \to 0$.
Fixing a decomposition $\cdop^n = \cdop \oplus \cdop^{n-1}$,
we obtain a section $\tilde s \in H^0(Q)$.
The composition with $\alpha$ gives a global section $s$ of $Hom(K,Q)$ the tangent bundle of $\Grass(n,k)$.
It is easy to see that the vanishing locus of this section is the disjoint union of $\Grass(n-1,k)$ and $\Grass(n-1,k-1)$.
Eventually, we obtain $\chi(\Grass(n,k))=\chi(\Grass(n-1,k))+\chi(\Grass(n-1,k-1))$.
The obvious fact, that $\Grass(n,n) \cong \Grass(n,0) \cong \Spec(\cdop)$ yields
$\chi(\Grass(n,k))=\binom{n}{k}$.
The decomposition could be used to define a $\cdop^*$-action on the Grassmannian
with fixed locus exactly $V(s)$.
Thus, the above can be seen as a special case of equivariant cohomology (see \cite{AB}).

\section{Local properties of reducible subschemes}\label{Slocal}
Let $V(s) = Z_1 \cup Z_2 \cup \ldots \cup Z_M$ be a decomposition of the vanishing locus into its irreducible components. We make the following

\neu{local1}
{\bf Assumption.}
{\em The $Z_i$ are locally complete intersections of codimension $d_i$.
Furthermore, the intersection $Z_{i_1} \cap Z_{i_2} \cap \ldots \cap Z_{i_k}$
is empty or a locally complete intersection of the expected codimension
$d_{i_1}+d_{i_2} + \ldots +d_{i_k}$.}

\neu{local2}
We want to compare the blow ups of $Z_1$ and $V(s)$ along $Z$.
Let us fix the notations.
We consider the local ring $A=\Ocal_{X,x}$ of a point $x$ in $X$.
We assume that the components $Z_1,Z_2,\ldots,Z_N$ of $V(s)$ pass through $x$.
By $Z$ we denote the component $Z_1$,
and by $Y$ we denote the union of $Z_2,\ldots,Z_N$.
We have the exact sequence
$$0 \to \Jcal_{V(s)} \to \Jcal_Z \to \Jcal_{Z,V(s)} \to 0 \,.$$

\neu{local3}{\bf Lemma.}
{\em The above exact sequence remains exact when twisted with $\Ocal_Z$.}

\proof We use the Koszul complex
$A^{\oplus \binom{d_1}{k}} \to A^{\oplus \binom{d_1}{k-1}}$ of $\Ocal_Z$
to compute the groups $\Tor_i^A( - , \Ocal_Z)$ (cf. \S 16 in \cite{Mat}).
Since the $d_1$ defining equations of $Z$ vanish on $Z$ and on $Z \cap Y$ we find
$$\Tor_i^A(\Ocal_Z , \Ocal_Z) = \Ocal_Z^{\oplus \binom{d_1}{i}} \qquad \qquad
\Tor_i^A( \Ocal_{Z \cap Y} , \Ocal_Z) =\Ocal_{Z \cap Y} ^{\oplus \binom{d_1}{i}}\, .$$
By our assumption $Z \cap Y$ is a locally complete intersection in $Y$.
Thus, $\Tor_i^A( \Ocal_Y,\Ocal_Z)=0$ for $i>0$.
Using these identities and the Mayer-Vietoris sequence
$$0 \to \Ocal_{Z \cup Y} \to \Ocal_Y \oplus \Ocal_Z \to \Ocal_{Z \cap Y} \to 0$$
we obtain
$\Tor_i^A( \Ocal_{Z \cup Y},\Ocal_Z)= \Jcal_{Z \cap Y,Z}^{\oplus \binom{d_1}{i}}$, for $i>0$,
and $\Ocal_{Z \cup Y} \otimes \Ocal_Z = \Ocal_Z$.
We conclude that $\Jcal_{Z \cup Y} \otimes \Ocal_Z \cong \Jcal_{Z \cap Y,Z}^{\oplus d_1}$.
We eventually obtain, that the torsion free $\Ocal_Z$-modul $\Jcal_{Z \cup Y} \otimes \Ocal_Z$ is a submodul of $J_Z/J^2_Z$,
because $\Tor_1^A(J_{Z,Z \cup Y},\Ocal_Z)$ is concentrated in $Z \cap Y$.
\qed

\neu{local4}
We conclude that the subalgebra $\oplus_{n \geq 0} \Jcal^n_{V(s)}$ of $\oplus_{n \geq 0} \Jcal^n_Z$ remains a subalgebra when restricted to $Z$.
Thus, the restriction of the blow up $\tilde X$ to $Z$ is a blow up of
the projective conormal bundle $\pdop(\Jcal_{Z_1}/\Jcal^2_{Z_1})$ of $Z_1$.

Let $D = D_1 \cup D_2 \cup \ldots \cup D_M$ be the decomposition of the exceptional divisor of the blow up $\tilde X \to X$ with $D_i$ dominating $Z_i$.
We have the following diagram
$$\xymatrix{& D_i \ar[d] \ar[r] \ar[dl]_-{\beta_i} & \tilde X \ar[d]^\pi \\
\pdop(\Jcal_{Z_i}/\Jcal^2_{Z_i}) \ar[r]^-{\pi_i} & Z_i \ar[r]^-{\iota} & X}$$

\neu{local5}{\bf Lemma.}{\em 
Let $c$ be a class in $A_m(X)$. Using the above notations we have}
$$\pi_*((c_1^k(\Ocal_{\tilde X}(D_i))\cap \pi^*(c))|_{D_i}) =
- \iota_*( (-1)^{d_1} s_{k+1-d_i}(\Jcal_{Z_i}/\Jcal^2_{Z_i}) \cap \iota^*c) \, .$$
\proof Since $\beta_i$ is a blow up and all the classes which appear in the formula are pull back classes from $\pdop (\Jcal_{Z_i}/\Jcal^2_{Z_i})$ we have (cf. Proposition 6.7 (b) in \cite{Ful}):
$$\begin{array}{rclr}
\pi_*((c_1^k(\Ocal_{\tilde X}(D_i)).\pi^*(c))|_{D_i})
&=& \iota_* \pi_{i*} (c_1^k(\Ocal_{\Jcal_{Z_i}/\Jcal^2_{Z_i}}(-1)) \cap \pi_i^* \iota^* c) \\
& = & (-1)^k \iota_* \pi_{i*} (c_1^k(\Ocal_{\Jcal_{Z_i}/\Jcal^2_{Z_i}}(1)) \cap \pi_i^* \iota^* c) \\
& = & (-1) ^{d_i-1} \iota_* (s_{k+1-d_i}(\Jcal_{Z_i}/\Jcal^2_{Z_i}) \cap \iota^* c) \, . & \quad \square
\end{array}$$

\section{The main formula}\label{Smain}
\neu{main}{\bf Theorem.}
{\em Let $E$ be a vector bundle on a equidimensional scheme $X$.
Let $s$ be a global section of $E$ whose vanishing locus $V(s)$ satisfies assumption \ref{Slocal}.\ref{local1}.
The evaluation of the top Chern class $c(E)$ at the fundamental cycle of $X$
can be computed along the components $Z_i$ of $V(s)$ and their intersections using the formula 
$$\int_Xc(E) = \sum_{ \{i_1,i_2,\ldots, i_k \} \subset \{ 1,2, , \ldots , M \} } -(-1)^k \int_{Z_{i_1} \cap Z_{i_2} \cap \ldots \cap Z_{i_k}} c(E|_{Z_{i_1} \cap Z_{i_2} \cap \ldots \cap Z_{i_k}}) . c(i_1,i_2,\ldots,i_k) \, ,$$
where the contribution $c(i_1,i_2,\ldots,i_k)$ of $Z_{i_1} \cap Z_{i_2} \cap \ldots \cap Z_{i_k}$ is given by

\tip
$c(i_1,i_2,\ldots,i_k) = 
\sum_{l_1 \geq 0, \, l_2 \geq 0, \ldots , \, l_k \geq 0} \binom{l_1+d_{i_1}+l_2+ d_{i_2} +\ldots + l_k+d_{i_k}}{l_1+d_{i_1},l_2+d_{i_2},\ldots,l_k+d_{i_k}}
s_{l_1}(N_{i_1}).s_{l_2}(N_{i_2})\ldots s_{l_k}(N_{i_k})$.

\tip
here $N_i$ denotes normal bundle of $Z_i$ restricted to $Z_{i_1} \cap Z_{i_2} \cap \ldots \cap Z_{i_k}$.
}

\proof
We start with the formula of Proposition \ref{Sreduction}.\ref{reduction2}.
Using the decomposition $D=D_1 \cup D_2 \cup \ldots \cup D_M$ we have $M$ contributions to $\int_Xc(E)$ by integrating the term over each individual divisor $D_i$.
Without loss of generality we set $i=1$.
By $\bar D_i$ we denote the first Chern class of $D_i$ restricted to $D_1$.
Then we have
$$\begin{array}{rcl}
\lefteqn{\int_{D_1} c(\pi^*E).c(\Ocal(D)|_{D_1})^{-1} = 
\int_{D_1} c(\pi^*E). \frac{1}{1+ (\bar D_1 + \bar D_2 + \ldots +\bar D_M) }}\\
\displaystyle & = & \displaystyle \int_{D_1} c(\pi^*E). \left(
\sum_{k \geq 0} (-\bar D_1 -\bar D_2 - \ldots -\bar D_M)^k
\right)\\
&=&\displaystyle \int_{D_1} c(\pi^*E). \left( \sum_{k_1 \geq 0, k_2 \geq 0, \ldots, k_M \geq 0} (-1)^{k_1+k_2 + \ldots + k_M} \binom{k_1+k_2 + \ldots + k_M}{k_1,k_2,\ldots,k_M} \bar D_1^{k_1}. \bar D_2^{k_2} \ldots \bar D_M^{k_M}
\right) \,.
\end{array}$$
For $2 \geq i \geq M$ we consider the exponent $k_i$ of $\bar D_i$.
If $k_i$ is positive, then we prefer to restrict to the intersection with $D_i$ which allows us to reduce the exponent $k_i$ by one.
On $D_1 \cap D_2 \cap \ldots \cap D_l$ we have a contribution of
$$\begin{array}{rcl}
\lefteqn{a:=\int_{D_1 \cap D_2 \cap \ldots \cap D_l} c(\pi^*E) \cdot}\\
& &\displaystyle {} \cdot \left(
\sum_{k_1 \geq 0, k_2 \geq 0, \ldots, k_l \geq 0}
(-1)^{l-1+k_1+k_2 + \ldots + k_l} 
\binom{l-1+k_1+k_2 + \ldots + k_l}{k_1,k_2+1,\ldots,k_l+1}
\bar D_1^{k_1}. \bar D_2^{k_2} \ldots \bar D_l^{k_l}
\right)
\end{array}$$

to the above sum.
Using lemma \ref{Slocal}.\ref{local5} $l$-times we obtain
$$\begin{array}{rcl}
\lefteqn{a= \int_{Z_1 \cap Z_2 \cap \ldots \cap Z_l} c(\iota^*E). \left(
\sum_{k_1 \geq 0, k_2 \geq 0, \ldots, k_l \geq 0}
(-1)^{d_1+d_2+\ldots d_l -1+k_1+k_2 + \ldots + k_l} \cdot \right. }\\
\\
& &
{} \cdot \left. \displaystyle \binom{l-1+k_1+k_2 + \ldots + k_l}{k_1,k_2+1,\ldots,k_l+1}
s_{k_1+1-d_1}(N_1\dual).s_{k_2+1-d_2}(N_2\dual) \ldots s_{k_l+1-d_l}(N_l\dual)
\right)
\end{array}$$
Applying $s_k(N)=(-1)^ks_k(N\dual)$, we obtain
$$\begin{array}{rl}
\lefteqn{a= \int_{Z_1 \cap Z_2 \cap \ldots \cap Z_l} c(\iota^*E). \left(
\sum_{k_1 \geq 0, k_2 \geq 0, \ldots, k_l \geq 0}
(-1)^{l -1} \cdot {} \right.}\\
& \displaystyle {} \cdot \left.
\binom{l-1+k_1+k_2 + \ldots + k_l}{k_1,k_2+1,\ldots,k_l+1}
s_{k_1+1-d_1}(N_1).s_{k_2+1-d_2}(N_2) \ldots s_{k_l+1-d_l}(N_l)
\right) \, .\end{array}$$
Since there are no negative Segre classes, we may change the summation to obtain:
$$\begin{array}{rl}
\lefteqn{a= \int_{Z_1 \cap Z_2 \cap \ldots \cap Z_l} c(\iota^*E). \left(
\sum_{k_1 \geq 0, k_2 \geq 0, \ldots, k_l \geq 0}
(-1)^{l -1} \cdot {} \right.}\\
& \displaystyle {} \cdot \left. \binom{-1+k_1+d_1+k_2+d_2 + \ldots + k_l+d_l}{k_1+d_1-1,k_2+d_2,\ldots,k_l+d_l}
s_{k_1}(N_1).s_{k_2}(N_2) \ldots s_{k_l}(N_l)
\right) \, .
\end{array}$$
Now the basic formula \ref{Snotations}.\ref{multinomial} for multinomial coefficients yields the theorem.
\qed

\neu{mainadd}{\bf Geometric interpretation of the terms $c(i_1,i_2,\ldots,i_k)$.}
First we remark that the $k$-th Segre class of a rank $r$ vector bundle $E$ could be defined by 
$$s_k(E)(\alpha) := (-1)^k \pi_*(c_1(\Ocal_E(1))^{r+k} \cap \pi^*(\alpha)) $$
where $\pi:\pdop(E \oplus \Ocal_X) \to X$ is the projective bundle associated to $E \oplus \Ocal_X$.
This definition yields $s_0(E)=1$. This fact allows the inversion of the total Segre class which gives the total Chern class of $E$.

This definition may be generalized to a finite set $E_1, \, E_2, \, \ldots , \, E_k$ of vector bundles:
We consider the
$\pdop^{r_1} \times \pdop^{r_2} \times \ldots \times \pdop^{r_k}$-fibre bundle
$$Y= \pdop(E_1 \oplus \Ocal_X) \times_X \pdop(E_2 \oplus \Ocal_X) \times_X  \ldots \times_X \pdop(E_k \oplus \Ocal_X)$$ over $X$, where $r_i=\rk(E_i)$.
On $Y$ we consider the relative ample bundle $\Ocal_Y(1)$ which is the product of the pull backs of the $\Ocal_{E_i}(1)$ to $Y$.
Let $r=r_1+r_2+\ldots + r_k$ be the relative dimension of the morphism $\pi:Y \to X$.
We define the $l$-th Segre class $s_l(E_1,E_2, \ldots ,  E_k)$ of the bundles 
$E_1, \, E_2, \, \ldots , \, E_k$ by 
$$s_l(E_1,E_2, \ldots ,  E_k)(\alpha) := (-1)^l \pi_*(c_1(\Ocal_E(1))^{r+l} \cap \pi^*(\alpha)) \, .$$
As usual, we define the total Segre class $s(E_1,E_2, \ldots ,  E_k)$ to be the sum $\sum_{l \geq 0} s_l(E_1,E_2, \ldots ,  E_k)$.

Using the notations of Theorem \ref{Smain}.\ref{main} we have that the contribution $c(i_1,i_2,\ldots,i_k)$ of the intersection
$Z_{i_1} \cap Z_{i_2} \cap \ldots \cap Z_{i_k}$ can be identified with the total Segre class of the restricted normal bundles:
$$c(i_1,i_2,\ldots,i_k) = s(N_{i_1},N_{i_2}, \ldots, N_{i_k}) \, .$$

\section{Counting rational curves on complete intersections}\label{Scounting}
\neu{counting1}
To compute the (virtual) number of curves of degree $k$ on a divisor $D$ of degree $d$ in $\pdop^n$ we consider the diagram
$$\xymatrix{\Mcal_0(\pdop^n,k) & \Ccal \ar[l]_-{q} \ar[r]^-{p} & \pdop^n}\, .$$
Here $\Mcal_0(\pdop^n,k)$ is the moduli space of stable mapping of rational curves to $\pdop^n$.
We have a 1-1 correspondence between sections of $\Ocal_{\pdop^n}(d)$ and sections of $q_*p^*\Ocal_{\pdop^n}(d)$.
If the rank of $q_*p^*\Ocal_{\pdop^n}(d)$ equals the dimension of the moduli space $\Mcal_0(\pdop^n,k)$, then its top Chern number is the virtual number of rational curves of degree $k$ on a divisor of degree $d$ in $\pdop^n$.

We propose the following: Take $D$ to be a divisor of $d$ hyperplanes which are in general position.
It follows that the corresponding section of $q_*p^*\Ocal_{\pdop^n}(d)$ vanishes on $d$ components which are isomorphic to $\Mcal_0(\pdop^{n-1},k)$.
Each of these components is smooth of codimension $k+1$ (see \cite{FP}).
Furthermore, the intersection of $l$ of these components is isomorphic to a 
$\Mcal_0(\pdop^{n-l},k)$ which is of codimension $l(k+1)$.
Thus, the conditions of assumption \ref{Slocal}.\ref{local1} are satisfied.

\neu{counting1a}
The bundle $N:=q_*p^*\Ocal_{\pdop^n}(1)$ restricted to $\Mcal_0(\pdop^{n-l},k)$ is the restriction of the normal bundle of one of the components (cf. 1.2 in \cite{hei}). Furthermore, all of these components have the codimension $k+1$.
Thus, the contribution $c(1,2,\ldots,l)$ can be simplified.
The following table shows some of these contributions.

$$\begin{array}{r|r|l}
k & l &  c(1,2,\ldots,l) \\
\hline
1 & 1 & 1+s_1(N)+s_2(N)+s_3(N) + \ldots = c(N)^{-1} \\
1 & 2 & 6 +20 s_1(N) + (20s_1^2(N)+30s_2(N)) +
(42s_3(N)+70s_1(N)s_2(N))+\ldots \\
1 & 3& 90 +630s_1(N) + (1680s_1^2(N)+1260s_2(N))+ {}\\&&\hfill {} + (2268s_3(N)+7560s_1(N)s_2(N)+1680s_1^3(N))+\ldots \\
2 & 1 &1+s_1(N)+s_2(N)+s_3(N) + \ldots = c(N)^{-1} \\
2 & 2 & 20 + 70 s_1(N) +(70s_1^2(N)+112s_2(N)) +(168s_3(N)+252s_1(N)s_3(N)) + \ldots \\
2 & 3 & 1680 +12600 s_1(N) + (27720s_2(N)+34650s_1^2(N)) + {} \\
&& \hfill {} + (55440s_3(N) + 166320s_1(N)s_2(N) + 34650s_1^3(N) +\ldots \\ 
\end{array}$$

\neu{counting2}{\bf Remark.}
The above considerations generalize to the case of complete intersections.
To consider a complete intersection of type $(d_1,d_2,\ldots,d_l)$ we have to consider the vector bundle $q_*p^*( \Ocal_{\pdop^n}(d_1) \oplus \Ocal_{\pdop^n}(d_2) \oplus \ldots \oplus \Ocal_{\pdop^n}(d_l) )$.
Moreover, we can apply our technique to complete intersections in arbitrary schemes provide we can degenerate these complete intersections to unions of convex schemes whose intersections are convex too.
Convexity is a condition which guarantees the smoothness of the moduli spaces 
under consideration (cf. \cite{FP}).

\neu{exam3}{\bf Lines on a cubic surface}
We consider the following situation:
$$\xymatrix{ \Grass(4,2) & \Ccal \ar[r]^-{p} \ar[l]_-{q} & \pdop^3}$$
where $\Ccal$ is the universal line over the Grassmannian scheme $\Grass(4,2)$ of lines in $\pdop^3$.
The vector bundle $E:=q_*p^* \Ocal_{\pdop^3}(3)$
is of rank 4 and the global sections of $E$ are in 1-1 correspondence with global sections of $\Ocal_{\pdop^3}(3)$.
Therefore $\int_{\Grass(4,2)}c(E)$ gives the number of lines on a cubic provided this number is finite.
It is not too hard to compute this number directly,
but we will do that using residual intersection theory with reducible cycles.

We choose a highly singular section of $\Ocal_{\pdop^3}(3)$ which consists of three planes which intersect in a point.
The corresponding section $s$ of $E$ vanishes in the three Grassmannians $\Grass(3,2)$ corresponding to the lines contained in each of the planes.
Each two of these Grassmannians intersect in a point.
To compute the contribution of each $\Grass(3,2)$ we have to regard the diagram
$$\xymatrix{\Grass(3,2) & \Ccal' \ar[r]^-{p} \ar[l]_-{q} & \pdop^2}.$$
Obviously $\Grass(3,2)$ is isomorphic to $\pdop^2$.
By $h$ we denote the hyperplane class.
The rank two bundle $N:=q_*p^*\Ocal_{\pdop^2}$ is isomorphic to $\Theta_{\pdop^2}(-1)$ and has the total Chern class $c(N)=1+h +h^2$.
The restriction of $E$ to $\Grass(3,2)$ is the symmetric product $\Sym^3N$.
Thus, we have $c(E|_{\Grass(3,2)})= 1+6h+21h^2$.
The normal bundle of $\Grass(3,2)$ in $\Grass(4,2)$ is $N$ (cf. 1.2 in \cite{hei}).
Hence the contribution of each $\Grass(3,2)$ is:
$$\int_{\Grass(3,2)} c(E|_{\Grass(3,2)}) \cdot c(N)^{-1} =
\int_{\pdop^2} (1+6h+21h^2)(1-h) =
\int_{\pdop^2} (1+5h+15h^2) =15 \, .$$
To compute the contribution of each of the intersection points we just have to
know that they are the intersection of two subschemes of codimension two.
By \ref{Smain}.\ref{main} the contribution of each point is $-\binom{4}{2}=-6$.
There are three $\Grass(3,2)$ and three intersection points which leads to the conclusion: There are 27 lines on a smooth cubic surface.

\section{Lines on a quintic threefold}\label{Sexam4}
\neu{exam41}
We consider the following universal family of lines in $\pdop^4$:
$$\xymatrix{ \Grass(5,2) & \Ccal \ar[r]^-{p} \ar[l]_-{q} & \pdop^4}\, .$$
As in Example \ref{Scounting}.\ref{exam3}  we see that the number of lines on a
general quintic is $\int_{\Grass(5,2)}c(E)$, where $E$ is the vector bundle $q_*p^*\Ocal_{\pdop^4}(5)$.
As before we reduce to calculations on Grassmannians of lower dimension  by taking the section of $E$ corresponding to the $\Ocal_{\pdop^4}(5)$ section with vanishing divisor 5 hyperplanes in general position.
Thus, we obtain three different contributions:

\begin{itemize}
\item 5 $\Grass(4,2)$-contributions of lines in one of these hyperplanes.
\item $\binom{5}{2}$ $\Grass(3,2)$-contributions of lines in intersections of two of these hyperplanes;
\item $\binom{5}{3}$ $\Grass(2,2)$-contributions of lines in intersections of three of these hyperplanes;
\end{itemize}
Before we start to compute these contributions we define the rank two vector bundle $N$ to be $q_*p^* \Ocal_{\pdop^4}(1)$.
Its fifth symmetric power is the vector bundle $E$.

\neu{exam42}{\bf The $\Grass(4,2)$-contribution.}
The restriction $Q$ of $N$ to $\Grass(4,2)$ is the tautological quotient.
Furthermore, it is the normal bundle to $\Grass(4,2)$ in $\Grass(5,2)$.
The Chern classes of $Q$ are called the special Schubert classes and denoted by $\sigma_i:=c_i(Q)$. We will need the intersection
numbers\footnote{To compute these intersection numbers we take the 
short exact sequence $0 \to K \to \Ocal_{\Grass(4,2)}^{\oplus 4} \to Q \to 0$.
Formal computation of $c(K)=c(Q)^{-1}$ gives
$c(K)=1-\sigma_1+(\sigma_1^2-\sigma_2) +(2\sigma_1\sigma_2-\sigma_1^3)+(\sigma_2^2-3\sigma_1^2\sigma_2+\sigma_1^4)$.
Since the rank of $K$ is two,
we obtain the two relations
$2\sigma_1\sigma_2-\sigma_1^3=0$,
and $\sigma_2^2-3\sigma_1^2\sigma_2+\sigma_1^4=0$.
Since $\Grass(4,2)$ is embedded as a quadric into $\pdop^5$ by the Pl\"ucker embedding, we have  $\int_{\Grass(4,2)}\sigma_1^4=2$.
The two relations and $\int_{\Grass(4,2)}\sigma_1^4=2$ yield the remaining two intersection numbers.}
$$\int_{\Grass(4,2)}\sigma_1^4=2 \qquad \int_{\Grass(4,2)}\sigma_1^2\sigma_2 = 1
\qquad \int_{\Grass(4,2)} \sigma_2^2 = 1 \, .$$
The restriction of $E$ to $\Grass(4,2)$ is $\Sym^5(Q)$.
Using the splitting principle we obtain
$$c(\Sym^5(Q))= 1+15\sigma_1+(85\sigma_1^2+35\sigma_2) +(225\sigma_1^3+350\sigma_1\sigma_2) +(274\sigma_1^4+1183\sigma_1^2\sigma_2+259\sigma_2^2) \, .  $$
Now we have anything at hand to compute the $\Grass(4,2)$-contribution:
$$\begin{array}{rl}
\int_{\Grass(4,2)} c(\Sym^5(Q)) \cdot c(N)^{-1}
&= \int_{\Grass(4,2)} c(\Sym^5(Q)) \cdot (1-\sigma_1+(\sigma_1^2-\sigma_2)) \\
& =\int_{\Grass(4,2)} (134 \sigma_1^4+783 \sigma_1^2\sigma_2 + 224\sigma_2^2 ) \\
& = 268 + 783+ 224 = 1275 \, .\\
\end{array}$$

\neu{exam43}{\bf The $\Grass(3,2)$-contribution.}
The restriction of $N$ to $\Grass(3,2)$ is again the tautological quotient bundle $Q$. Using the identification $\Grass(3,2) \cong \pdop^2$ we find $c(N)=1+h +h^2$ where $h$ is the hyperplane class in $\pdop^2$.
Thus, we have $s_1(N) =-h$, and $s_2(N)=0$. 
$N$ is the restriction of the normal bundle of the two $\Grass(4,2)$ which cut out the $\Grass(3,2)$.
The contribution $c(1,2)$ of Theorem \ref{Smain}.\ref{main} can now easily be computed using the formulas of \ref{Scounting}.\ref{counting1a}.
$$c(1,2) =  6 + 20 s_1(N)+20s_1^2(N) + 30 s_2(N) = 6- 20 h +20 h^2 \, .$$
Again we have that the restriction of $E$ to $\Grass(3,2)$ is $\Sym^5(N)$.
The splitting principle yields
$$c(\Sym^5(N)) = 1+15c_1(N) +(85c_1^2(N)+35c_2(N)) = 1+15h + 120h^2 \, .$$
Using Theorem \ref{Smain}.\ref{main} we obtain the $\Grass(3,2)$-contribution to be
$$\int_{\Grass(3,2)} c(E|_{\Grass(3,2)}) \cdot c(1,2) = \int_{\pdop^2} (1+15h + 120h^2)(6- 20 h +20 h^2) = 440 \, .$$

\neu{exam44}{\bf The $\Grass(2,2)$-contribution.}
Since $\Grass(2,2)$ is a point we just have to take the starting term of $c(1,2,3)$ to get its contribution.
Thus, the table of \ref{Scounting}.\ref{counting1a} tells us $c(1,2,3)=90$.

\neu{exam45}{\bf Drawing the conclusion.}
Summing up we obtain, that there are
$$5 \cdot 1275 -\binom{5}{2}\cdot 440+\binom{5}{3}\cdot 90 = 6375-4400+900=2875$$
lines on a general quintic threefold.

\vspace{5em}
{\small
Georg Hein\\
Humboldt-Universit\"at zu Berlin\\
Institut f\"ur Mathematik\\
Unter den Linden 6\\
10099 Berlin, Germany\\
{}\\
hein@mathematik.hu-berlin.de}
\end{document}